\documentclass[a4paper,12pt]{article}
\usepackage{amsmath}
\usepackage{amsfonts}
\usepackage{amssymb}

\advance \topmargin by -\headheight
\advance \topmargin by -\headsep
\newtheorem{remark}{Remark}[section]
\newtheorem{theorem}{Theorem} [section]
\newtheorem{lemma}{Lemma}[section]

\newcommand{\be}{\begin{equation}} 
\newcommand{\ee}{\end{equation}}
\newcommand{\bea}{\begin{eqnarray}} 
\newcommand{\eea}{\end{eqnarray}}
\newcommand{\bean}{\begin{eqnarray*}} 
\newcommand{\eean}{\end{eqnarray*}}
\newcommand\R{{\mathbb R}}
\newcommand\N{{\mathbb N}}
\newcommand\Z{{\mathbb Z}}

\def\Proof{{\smallskip\noindent{\em Proof. }}}
\def\endProof{{\hfill$\Box$\medskip\noindent}}

\begin{document}

\title{Global existence versus blow up \\for some models 
        of interacting particles}  
\author{Piotr Biler$^1$ and Lorenzo Brandolese$^2$\\
\small $^1$ Instytut Matematyczny, Uniwersytet Wroc{\l}awski, \\
\small pl. Grunwaldzki 2/4, 50--384 Wroc{\l}aw, POLAND\\
\small{\tt Piotr.Biler@math.uni.wroc.pl}\\
\small $^2$ Institut Camille Jordan,  Universit\'e Claude Bernard -- Lyon 1,\\
\small  21 avenue Claude Bernard,  69622 Villeurbanne Cedex, FRANCE \\
\small {\tt brandolese@math.univ-lyon1.fr}}

\maketitle

\begin{abstract}
\noindent
We study the global existence and space-time asymptotics of  solutions 
for a  class of nonlocal parabolic semilinear equations.
Our models include the Nernst--Planck and the Debye--H\"uckel drift-diffusion systems 
as well as parabolic-elliptic systems of chemotaxis.
In the case of a~model of self-gravitating particles, 
we also give a result on the finite time blow up
of solutions with localized and oscillating com\-plex-valued initial data, 
using a method by S. Montgomery-Smith.
\end{abstract}

\section{Introduction}

In this paper we are concerned with semilinear parabolic systems of the form
\begin{eqnarray}
\label{general syst}
\partial_t u_j &=&\Delta u_j +\nabla\cdot\left(
\textstyle\sum_{h,k=1}^m c_{j,h,k}\,u_h(\nabla E_d*u_k)\right), \quad j=1,\ldots,m,\\
u(0)(x)&=&u_0(x).\nonumber
\end{eqnarray}
Here the unknown is the vector field $u=(u_1,\ldots,u_m)$, defined on the whole space $\R^d$ (with $m\ge1$ and $d\ge2$),
and $c_{j,h,k}\in L^\infty(\R^d)$, $j,h,k=1,\ldots,m$, are given coefficients.
Moreover, $E_d$ denotes the fundamental solution of the Laplacian in $\R^d$.

Systems of the form \eqref{general syst} arise e.g. from plasma, semiconductors and electrolytes theories, 
                  biology (modelling of chemotaxis phenomena) and statistical mechanics. 
The basic example for us is the model for gravitating particles: in this case $m=1$, $d\ge 2$, and 
the governing equations are usually written as
\begin{eqnarray}
\label{gravitating}
\partial_t u=\Delta u+\nabla\cdot(u\nabla\varphi), \qquad  
\Delta\varphi=u.
\end{eqnarray}
Here $u=u(x,t)$ is the density of the particles
and  $\varphi$ is the self-consistent gravitational potential generated by $u$. 
Related systems appear also in the theory of chemotaxis, see e.g. \cite{CPZ},  \cite{BDP}, \cite{B-AMSA}. 
We do not require that $u\ge0$ in our study, which is, however, relevant in physical applications;  
we even admit complex-valued solutions.
In this case the coefficients $(c_{j,h,k})$ are constant and equal to~$1$.
Another important example is provided by the Debye system, in which
the first equation of \eqref{gravitating} is replaced with
\be
\partial_t u=\Delta u-\nabla\cdot(u\nabla\varphi).\label{Debye}
\ee

A more general model,  still belonging to the class \eqref{general syst},
is the drift-diffusion system
\begin{eqnarray}
\label{drift diff}
\partial_t v&=&\Delta v-(\nabla \cdot(v\nabla \phi)),\nonumber\\
\partial_t w&=&\Delta w+(\nabla \cdot(w\nabla\phi)),\\
\Delta\phi&=&v-w.\nonumber
\end{eqnarray}
In the theory  by W.~Nernst and M.~Planck, $v$ and $w$  represent the density of positively  
and negatively charged particles, respectively.

A lot is known about the existence and the nonexistence of real-valued solutions of these models,
see e.g. \cite{BHN}, \cite{BCGK04}, \cite{B-SM}, \cite{Bil95}, and the references therein.
For instance, if $d=1$, then the considered models have global in time solutions. 
If $d\ge 2$, the Debye system \eqref{Debye} and more general \eqref{drift diff} 
have global in time solutions and their asymptotics is described by suitable self-similar 
solutions, \cite{BD00} and \cite{MO}. 
These may be interpreted as a complete diffusion of charges to infinity due to repulsive interactions. 

On the other hand, models describing 
either chemotaxis or gravitational interaction in $d\ge 2$ dimensions feature concentration 
phenomena which may eventually lead to a collapse of solutions. These phenomena manifest by 
the formation of singularities of solutions like weak convergence either to Dirac point masses or 
to unbounded functions $\sim |x|^{-2}$. 

One purpose of this paper is to show that a different kind of finite time blow up can occur
for solutions of \eqref{gravitating} (and for a few other particular cases of \eqref{general  syst}).
In particular, we will show that also 
\emph{nonpositive (in fact: complex-valued) and oscillating} solutions can blow up.
Our second purpose is to
give a global existence result for ``small'' solutions of \eqref{general syst}.
Such result will provide us with some decay profiles in space-time of solutions.

The global existence result for ``well localized'' solutions can be stated as follows
(see section \ref{section global} for a more general, and more precise, statement).

\begin{theorem}
\label{theorem1 intro}
Let $d\ge3$.
There exists $\eta>0$ such that if 
\begin{equation}
\label{small1 intro}
|u_0(x)|\le \frac{\eta}{(1+|x|)^2},
\end{equation}
then there exists $C\ge0$ and a  unique solution $u$ of~\eqref{general syst}
such that, for all $x\in\R^d$ and $t\ge0$,
\begin{equation}
|u(x,t)|\le {\frac{C}{(1+|x|)^2}}, \quad\hbox{ and }\quad
|u(x,t)|\le \frac{C}{1+t}.
\end{equation}
Moreover, if all the coefficients $c_{j,h,k}(x)$ are constant in $\R^d$, then the smallness
assumption~\eqref{small1 intro}, can be replaced by the weaker, scale invariant, condition
\begin{equation}
\label{small2 intro}
\hbox{\rm ess\,}\sup_{x\in \R^d} |x|^2|u_0(x)|\le \eta.
\end{equation}
\end{theorem}

It would be possible to establish similar decay profiles in space-time for the solution,
with a spatial decay rate  \emph{larger} than two.
In this case the decay rate as $t\to\infty$ is also increased, up to one-half of the decay
rate as $x\to\infty$ (or up to the rate $d/2$ if the space decay rate is larger than $d$).
The exponent two is however the most interesting case, since it
plays a special role in these models, for scaling reasons.
For example, it corresponds to the expected decay rate of self-similar solutions,
see  e.g. \cite{BCGK04}.
It is also the decay of the well-known Chandrasekhar solution,
$\tilde u(x,t)=2(d-2)|x|^{-2}$, which is a  stationary solution for \eqref{gravitating}
for $d\ge3$.
 
\medskip

We now state our result on the blow up of solutions.

\begin{theorem}
\label{theorem2 intro}
There exists $u_0\in {\cal S}_0(\R^d)$ 
(the space of functions belonging to the Schwartz class, with vanishing moments
of all order),
such that
the corresponding solution $u$ of~\eqref{gravitating} 
blows up in finite time: there exists $t^*>0$
such that $u(t^*)\not\in \dot B^{s,q}_p(\R^d)$
for all $s\in\R$, $1\le p,q\le\infty$.
 \end{theorem}

We will restate this theorem in a more precise way in section~\ref{section blow up}.
Therein, we will also recall the definition of the Besov norm $\|\cdot\|_{\dot B^{s,q}_p(\R^d)}$.
Here we only observe that this theorem tells us,  in particular, that
$\|u(t)\|_{L^p}$ blows up  for all $1\le p\le\infty$.
The proof of Theorem~\ref{theorem2 intro}
consists in proving suitable lower bound estimates for the Fourier transform $\widehat u$.
We shall derive such estimates using an idea of Montgomery-Smith \cite{MontS01}.

Our explosing solution $u$ of~Theorem~\ref{theorem2 intro} is in fact complex-valued 
since its Fourier transform enjoys some \emph {positivity and nonsymmetry} properties.
Of course, one  can rewrite the scalar equation~\eqref{gravitating} for the real and imaginary part of~$u$.
This yields  a blow up result for a real system of the form~\eqref{general syst},
which is formally close to~\eqref{drift diff}. 

\paragraph{Notations.}
In chains of inequalities, all the constants will be denoted by $C$ 
even if they vary from line to line.
We will simply write $\int$ instead of $\int_{\R^d}$.

\section{Global existence for the general model}
\label{section global}
The proof of Theorem~\ref{theorem1 intro}
relies on  size estimates of the kernel $\nabla E_d$.
We have of course
\begin{equation*}
|\nabla E_d(x)|\le \frac{C}{|x|^{d-1}}.
\end{equation*}
Our method applies also if we replace $\nabla E_d$ with 
any  kernel $K$ such that
$K$ is a~measurable function in $\R^d$, and
 \begin{equation}
 \label{decay K}
 |K(x)|\le C|x|^{-d-1+\alpha}, \qquad 1<\alpha<d.
 \end{equation}
In what follows, we will consider this more general situation.
For the applications that we have in mind, $\alpha=2$, and this explains
the restriction $d\ge3$ in Theorem~\ref{theorem1 intro}.
The two-dimensional case is often a special case
in these models, e.g. the Keller--Segel parabolic-elliptic model of chemotaxis.

Since we would obtain the same bounds for all the components of $u$,
in the remaing part of the section we can assume
that $u$ is scalar and $c_{j,h,k}=c$. Moreover, without loss of generality in the
proof below, can assume that $c$ is constant, essentially because
multipling $K$ with a $L^\infty$ function does not affect~\eqref{decay K}.

With this simplification, 
the models discussed above can be written in the following integral form
 \begin{equation}
 \label{IE}
 u(t)=e^{t\Delta}u_0+\int_0^t G(t-s)*\bigl(u(K*u)\bigr)(s)\,ds,
 \end{equation}
 where 
 $G$ behaves like a first order derivative of the Gaussian heat kernel, namely
 \begin{eqnarray}  
 \label{cond G}
&&G(x,t)=t^{-(d+1)/2}\Psi(x/\sqrt t),\quad\hbox{ with }\quad \Psi\in {\cal S}(\R^d),\\
\label{cond2 G}
&& e^{t\Delta/2}G\left(\frac{t}{2}-s\right)=G(t-s). 
 \end{eqnarray}

To state our result in a precise way, we introduce a few useful spaces.
 For $\theta\ge 0$ we define by $L^\infty_\theta$ the space of measurable functions
$f$ on $\R^d$, such that $(1+|\cdot|)^\theta f\in L^\infty(\R^d)$. 
 Let ${\cal E}_\theta$ the space of all measurable functions $f=f(x,t)$ in $\R^d\times\R^+$,
 such that
 \begin{eqnarray}
&& \hbox{ess}\sup_{\!\!\!\!\!\!\!\!\!\!\!\! x\in\R^d,\;t\ge0} (1+|x|)^\theta|f(x,t)|<\infty,\\
&& \hbox{ess}\sup_{\!\!\!\!\!\!\!\!\!\!\!\! x\in\R^d,\;t\ge0} (1+t)^{\theta/2}|f(x,t)|<\infty,
 \end{eqnarray}
 and
 \begin{equation} 
 \label{continuity}
 f\in {\cal C}((0,\infty), L^\infty_\theta).
 \end{equation}
The space ${\cal E}_\theta$ is equipped with its natural norm.
  
 \begin{theorem}
 \label{theorem1}
 Let $1< \alpha<d$ and $K$ such that \eqref{decay K} holds.
 Assume that $u_0\in L^\infty_{\alpha}$.
 Then we can find $\eta>0$ such that if
 \begin{equation}  
 \label{small1}
  \|u_0\|_{L^\infty_{\alpha}}<\eta,
 \end{equation}
 then there exists a unique mild solution $u$ of \eqref{IE},
 such that $u\in {\cal E}_\alpha$ and
  $u(t)\stackrel{{\cal D}'}{\longrightarrow} u_0$ as $t\to0$.
 Moreover,
 if we assume, in addition, that $K$ is homogeneous of degree $-d-1+\alpha$,
 then \eqref{small1} can be replaced by the  weaker, scale invariant condition:
\begin{equation}
\label{small2}
\hbox{\rm ess\,}\sup_{x\in\R^d} |x|^\alpha |u_0(x)|<\eta.
\end{equation}
\end{theorem}

The proof relies on the following simple lemma

\begin{lemma}
\label{lemma1}
Let $K$ satisfy the assumption  \eqref{decay K}.
\begin{enumerate}
\item 
If $f\in L^\infty_\alpha$, then $K*f\in L^\infty_1$.
\item
If $\phi\in {\cal E}_\alpha$, then $K*\phi\in {\cal E}_1$.
\end{enumerate}
\end{lemma}

\Proof
Using the duality and the interpolation of Lorentz spaces, we get
$$
     \|K*f\|_{L^\infty}\le \|K\|_{L^{d/(d+1-\alpha),\infty}} \|f\|_{L^{d/(\alpha-1),1}}
\le C\|f\|_{L^{d/\alpha,\infty}}^{1-1/\alpha}   \|f\|_{L^\infty}^{1/\alpha}.
$$
Thus,
$$ \|K*f\|_{L^\infty}\le C \|f\|_{L^\infty_\alpha}.$$
In particular, we may assume $|x|\ge1$.
Note that
$$ |K*f(x)|\le C\int |x-y|^{-d-1+\alpha}|f(y)|\,dy=I_1+I_2+I_3,$$
where 
$I_1\equiv\int_{|y|\le \frac{|x|}{2}}\dots\,dy$, 
$I_2\equiv\int_{\frac{|x|}{2} \le |y|\le \frac{3}{2}|x|}\dots\,dy$ and
$I_3\equiv\int_{|y|\ge \frac{3}{2}|x|} \dots\,dy$.
One easily checks that these three integrals are bounded by $C|x|^{-1}$.
The first part of the lemma follows.

On the other hand, by the above inequality,
$$
\|K*\phi(t)\|_{L^\infty}\le C \|\phi\|_{L^\infty(\R^+,L^{d/\alpha,\infty})}^{1-1/\alpha}   \|\phi(t)\|_{L^\infty}^{1/\alpha}
\le C(1+t)^{-1/2}\|\phi\|_{E_\alpha}.
$$
Combining this with the first part of the lemma, applied to $\phi(t)$, yields the result.
\endProof

\bigskip

For $u\in {\cal E}_\theta$, the nonlinear term $u(K*u)$ belongs, by Lemma~\ref{lemma1}, to
${\cal E}_{\theta+1}$.
Then it is natural to study the behavior of the linear operator
\begin{equation}
\label{def L}
L(w)(t)=\int_0^t G(t-s)*w(s)\,ds 
\end{equation}
in such a space.

\begin{lemma}
\label{lemma2}
Let $1<\alpha<d$ and $w\in {\cal E}_{\alpha+1}$.
Then $L(w)\in {\cal E}_\alpha$.
\end{lemma}

\Proof 
We will use repeatedly the property 
$$\|G(t-s)\|_{L^1}=C(t-s)^{-1/2},$$
which is a consequence of~\eqref{cond G}.
A few estimates below bear some relations with those of Miyakawa \cite{Mi00},
yielding space-time decay results for the Navier--Stokes equations.
We start observing that $L(w)\in L^\infty(\R^d\times\R)$.
Indeed,
\begin{eqnarray*}   
\| L(w)(t)\|_{L^\infty} &\le& \int_0^t \|G(t-s)\|_{L^1}\|w(s)\|_{L^\infty}\,ds\\
   & \le& C\|w\|_{{\cal E}_{\alpha+1}} \int_0^t (t-s)^{-1/2}s^{1/2}\,ds
    \le C\|w\|_{{\cal E}_{\alpha+1}}.
\end{eqnarray*}
Then we can assume in the following that $|x|\ge1$ and $t\ge1$.

We can write $L(w)=I_1+I_2$, where,
$$ I_1\equiv \int_0^t\!\!\int_{|y|\le |x|/2} G(x-y)w(y,s)\,dy\,ds$$
and
$$ I_2\equiv \int_0^t\!\!\int_{|y|\ge |x|/2} G(x-y)w(y,s)\,dy\,ds.$$
Now,
\begin{eqnarray*}
|I_1(x,t)| &\le & |x|^{-d}\int_0^t\!\!\int_{|y|\le |x|/2}(t-s)^{-1/2}(1+|y|)^{-\alpha}(1+s)^{-1/2}\,dy\,ds\\
  &\le&
  C|x|^{-d}\int_{|y|\le |x|/2} |y|^{-\alpha}\,dy\le C |x|^{-\alpha}.
\end{eqnarray*}
On the other hand,
\begin{equation*}
|I_2(x,t)| \le  C|x|^{-\alpha}\int_0^t \|G(t-s)\|_{L^1} \,s^{-1/2}\,ds\le
  C|x|^{-\alpha}.
\end{equation*}
Thus, $|L(w)(x,t)|\le C(1+|x|)^{-\alpha}\|w\|_{{\cal E}_{\alpha+1}}$
and, in particular,
$$ \|L(w)\|_{L^\infty((0,\infty),L^{d/\alpha,\infty})} \le C\|w\|_{{\cal E}_{\alpha+1}}.$$ 

To obtain a decay estimate as $t\to\infty$, we recall~\eqref{cond2 G} and
write
$$L(w)(t)=e^{t\Delta/2}L(w)(t/2)+\int_{t/2}^t G(t-s)*w\,ds\equiv J_1+J_2.$$
By duality (we denote here by $g_t$ the Gaussian kernel),
\begin{eqnarray*}
\|J_1\|_{L^\infty} &\le&
  C\|g_{t/2}\|_{L^{d/(d-\alpha),1}}\|L(w)(t/2)\|_{L^{d/\alpha,\infty}}, \\
&\le& Ct^{-\alpha/2}\|w\|_{{\cal E}_{\alpha+1}}.
\end{eqnarray*}
Moreover,
$$
\|J_2(t)\|\le Ct^{-(\alpha+1)/2}\|w\|_{{\cal E}_{\alpha+1}}\int_{t/2}^t \|G(t-s)\|_{L^1}\,ds
  \le Ct^{-\alpha/2}\|w\|_{{\cal E}_{\alpha+1}} .
  $$
The decay estimates in space-time for $L(w)$ then follow.
The continuity with respect to $t$ being straightforward,
the proof of Lemma~\ref{lemma2} is finished.
\endProof

\bigskip

By Lemma~\ref{lemma2}, the bilinear operator 
\begin{equation} 
\label{bilinear}
B(u,v)=\int_0^t G(t-s)*\bigl(u(K*v)\bigr)(s)\,ds
\end{equation}
is continuous from ${\cal E}_{\alpha}\times {\cal E}_{\alpha}$ to ${\cal E}_{\alpha}$.
Note that our last lemma also implies that $\|u(t)-e^{t\Delta}u_0\|_{L^\infty}\le C\sqrt t$,
so that $u(t)\to u_0$ a.e. and in the distributional sense.
The existence (and the uniqueness) of the solution of \eqref{IE}, under the assumption~\eqref{small1}
now follows by a standard argument, i.e. 
the application of the contraction mapping theorem.

\medskip
In order to finish the proof of Theorem~\ref{theorem1}
it only remains to show that the smallness assumption~\eqref{small1}
can be relaxed, when the kernel $K$ is a homogeneous function.
Consider the rescaling
\begin{equation}
\label{rescaling}
u_\lambda(x,t)=\lambda^{\alpha} u(\lambda x,\lambda^2 t).
\end{equation}
A direct computation shows that, if $K$ is homogeneous of degree $-d-1+\alpha$,
and $u$ is a solution of~\eqref{IE}, then $u_\lambda$ is a solution of~\eqref{IE} as well.
Now let $\eta>0$ be the constant obtained in the first part of Theorem~\ref{theorem1}.
Assume that the datum $u_0$ is such that~\eqref{small2} holds. Then we can choose 
a $\tilde \lambda>0$ such that
$$ \hbox{ess\,}\sup_{x\in\R^d}\tilde\lambda^\alpha(1+|x|)^\alpha|u_0(\tilde\lambda x)|<\eta.$$
We can apply the first part of Theorem~\ref{theorem1} to the initial datum
$\tilde\lambda^\alpha u_0(\tilde\lambda\,\,\cdot)$.
If we denote by $\tilde u$ the corresponding solution, we see that $\tilde u_{\tilde \lambda^{-1}}$
is the  solution of \eqref{IE} starting from $u_0$.
Theorem~\ref{theorem1} is now established.
\endProof 

\bigskip

\begin{remark}
\label{local solutions}
\begin{rm}
With minor modifications of  the decay exponents in the above proof,
one sees that, for any finite $T>0$
 the bilinear operator~\eqref{bilinear}
is bicontinuous in the space ${\cal C}([0,T], L^\infty_\theta)$,
for all $\theta\ge0$.
The contraction mapping theorem guarantees
that, if $u_0\in L^\infty_\theta$, $\theta\ge 0$ (with arbitrary norm) and $T>0$ is small enough,
then there exists a unique solution $u\in {\cal C}((0,T], L^\infty_\theta)$,
such that $u(t)\to u_0$ in the weak sense; 
we will write $u\in{\cal C}_w([0,T],L^\infty_\theta)$ to express these properties.
\end{rm}
\end{remark}

\section{Blow up for the model of gravitating particles}
\label{section blow up}

In this section we show that there exist solutions of~\eqref{gravitating},
with initial data $u_0$ in the Schwartz class,  and such that
$\int x^\alpha u_0(x)\,dx=0$ for all $\alpha\in \N^d$, which blow up in finite time.
Here we adopt a quite general definition of solution:
we ask that the Fourier transform $\widehat u(\cdot, t)$, also denoted $\widehat u_t$,
satisfies for a.e. $\xi\in\R^d$
and all $t\in[0,T]$, $0<T\le \infty$, the integral equation 
\begin{equation}
\label{fourier gravitating}
\widehat u_t(\xi)=e^{(s-t)|\xi|^2}\widehat u_0(\xi)
 +\frac{1}{(2\pi)^d}
 \int_0^t e^{(s-t)|\xi|^2}  i\xi \cdot  \biggl ( \widehat u_s(\xi) *\frac{ i\xi }{-|\xi|^2}\widehat u_s(\xi)\biggr)\,ds.
\end{equation}
The definition of the Fourier transform for integrable functions that we adopt is
$\widehat u_t(\xi)=\int u(x,t)e^{-i\xi\cdot x}\,dx$.

There are several ways to give a sense to the above integral and ensure the validity of~\eqref{fourier gravitating}.
An obvious way, is to consider the (local) solutions obtained in the setting of 
Remark~\ref{local solutions}, with $\theta>d$.
But the above equality is true in  more general settings.
For example, it holds
for the solutions $u\in {\cal C}_w([0,T]\colon {\cal PM}^{d-2})$,
(with $0<T\le \infty$ and $d\ge3$) constructed
in \cite{BCGK04}, where ${\cal PM}^{a}$ is the space of pseudomeasures 
$${\cal PM}^a=\{v\in{\cal S}'(\R^d)\colon \widehat v\in L^1_{\rm loc}(\R^d),\,
    \|v\|_{{\cal PM}^a}\equiv\hbox{ess}\!\sup_{\xi\in\R^d}|\xi|^a|\widehat v(\xi)|<\infty\}.$$
As pointed out in \cite{BCGK04}, a distributional solution of the Cauchy problem for \eqref{gravitating},
does also satisfy~\eqref{fourier gravitating}.

In this section we will consider initial data with nonnegative Fourier transform.
Under this condition, one immediately checks that the iteration scheme
yielding a solution in ${\cal C}_w([0,T],{\cal PM}^{d-2})$, or in ${\cal C}_w([0,T],L^\infty_\theta)$,
converges in the subset of functions $u$ such that $\widehat u(\xi,t)\ge0$, for
all $t\in[0,T]$ and a.e. $\xi\in\R^d$.
The crucial fact  that will lead to the blow up is the following:

\begin{lemma}\label{H}
Let
$$H_j(\widehat u)(\xi,t)\equiv\int_0^t\!\!\int e^{(s-t)|\xi|^2}\xi_j\frac{\eta_j}{ |\eta|^2}
    \widehat u_s(\xi-\eta)\widehat u_s(\eta)\,d\eta\,ds, 
    \qquad j=1,\ldots,d.$$
Then, if $0\le \widehat u\le \widehat v$, $\hbox{supp\,}\widehat u$ and 
$\hbox{supp\,}\widehat v$ are contained in $\{\xi\in\R^d\colon \xi_\ell\ge0,\, \ell=1,\ldots,d\}$,
then $0\le H_j(\widehat u)\le H_j(\widehat v)$, 
and their supports are still contained in $\{\xi\in\R^d\colon \xi_\ell\ge0,\, \ell=1,\ldots,d\}$. 
\end{lemma}

This simple  observation  allows us to adapt to our situation the argument
introduced by Montgomery-Smith for the ``cheap'' Navier--Stokes equations, 
see \cite{MontS01}.

Let us first recall the definition of the Besov norm $\|\cdot \|_{\dot B^{a,\infty}_\infty}$,
$a\in\R$.
We  consider a function $\psi\in{\cal S}(\R^d)$, such that
$\widehat \psi\ge0$ in $\R^d$, $\widehat \psi(\xi)\ge1$ for $\frac{1}{2}\le |\xi|\le 1$, 
$\widehat \psi(\xi)=0$ for $|\xi|\le\frac{1}{3}$ or $|\xi|\ge \frac{4}{3}$.
Then, for a distribution $f$, we define
\begin{equation}
\label{besov norm}
\|f\|_{\dot B^{a,\infty}_\infty}=\sup_{k\in\Z}2^{(a+d)k}\|\psi(\cdot \,2^k)*f\|_{L^\infty}.
\end{equation}

\begin{remark}
\label{other norm}
\begin{rm}
It is well-known that any Besov space $\dot B^{s,q}_p(\R^d)$,
as well as any Triebel--Lizorkin space $\dot F^{s,q}_p(\R^d)$
(so in particular the $L^p$-spaces, which are identified to $\dot F^{0,2}_p$),
 with 
$s\in\R$, $1\le p,\, q\le\infty$, are embedded in $\dot B^{s-{d/p},\infty}_\infty(\R^d)$.
It is then sufficient to show that $u(t)$ blows up in the $\dot B^{a,\infty}_\infty$ norm,
for all $a\in\R$, to deduce that all Besov and Triebel--Lizorkin norms of $u$ must blow up.
To be more precise, $L^1$ is not a Triebel--Lizorkin space, but
 we  will see that $\widehat u_{t^*}$ becomes unbounded for a finite $t^*$, hence
$\|u(t)\|_{L^1}$  does blow up.

A similar remark applies to pseudomeasure norms, since
$ {\cal PM}^a$ is continuously embedded in $\dot B^{a-d,\infty}_\infty(\R^d)$.
\end{rm}
\end{remark}

\begin{theorem}
\label{theorem2}
Let  $w_0\in{\cal S}(\R^d)$, such that $\widehat w_0$ is nonnegative and supported in the ball
$B_{1/4}(3e_1/4)$, where $e_1$ is the unit vector, and
$\|\widehat w_0\|_{L^1}=1$.
Let $A>2^{4/3}(2\pi)^d$ and $u_0=Aw_0$ (so in particular $u_0\in{\cal S}_0(\R^d)$).
Assume also that 
$u(\cdot, t)$ is a tempered distribution such that for all $t\ge0$, $\widehat u_t\ge0$ and~\eqref{fourier gravitating}
holds for a.e. $\xi\in \R^d$.
Then,  for all  $a\in\R$,
\begin{equation}
\|u(\cdot,t^*)\|_{\dot B^{a,\infty}_\infty}=\infty, \qquad\hbox{where\ \ $t^*=\log(2^{1/3})$}.
\end{equation}
\end{theorem}

\Proof
Set
$t_0=0$, $t_k=\log 2\left(\sum_{j=1}^k 2^{-2j}\right)$ and  $w_k=w_0^{2^k}$.
We set also
$$ \alpha_k(t)=2^{2k+7-2^k}e^{-2^kt}{\bf 1}_{t\ge t_k} \qquad
    (k\in\N),$$
and claim that, for $k=0,1,\ldots$,
\begin{equation}
\label{induction}
\widehat u_t(\xi)\ge A^{2^k}\alpha_k(t)\widehat w_k(\xi).
\end{equation}
This is seen by induction. For $k=0$ the claim follows from Lemma \ref{H}: 
$$\widehat u_t(\xi)\ge Ae^{-t|\xi|^2}\widehat w_0(\xi)\ge Ae^{-t}\widehat w_0(\xi),
\qquad t\ge0.$$ 
Now assume that~\eqref{induction}
holds for $k-1$.
Set
$$E_k=\{\xi\in\R^d\colon 2^{k-1}\le \xi_1\le |\xi|\le 2^k\}, \qquad k=0,1\ldots$$
Note that $\widehat w_k= (2\pi)^{-d}\widehat w_{k-1}*\widehat w_{k-1}$.
Thus, $\hbox{supp\,}\widehat w_k\subset E_k$.

But, for a.e. $\xi\in E_k$, estimating from  below by zero all the terms on
the right hand side of~\eqref{fourier gravitating},
except for the first term obtained after computing the scalar product, we get 
\begin{eqnarray*}
\widehat u_t(\xi) &\ge &
\frac{1}{(2\pi)^d}\int_0^t\!\!\int e^{(s-t)|\xi|^2} \frac{\xi_1\eta_1}{|\eta|^2} 
\widehat u_s(\xi-\eta)\widehat u_s(\eta)
\,d\eta\,ds\\
&\ge&
\int_{t_{k-1}}^t\!\!\int_{\eta\in E_{k-1}} e^{(s-t)|\xi|^2} \frac{\xi_1\eta_1}{|\eta|^2} 
 \bigl(A^{2^{k-1}}\alpha_{k-1}(s)\bigr)^2 
  \frac{\widehat w_{k-1}(\xi-\eta) \widehat w_{k-1}(\eta)}{(2\pi)^d }\,d\eta\,ds\\
 &\ge&
 A^{2^{k}}  2^{4k+8- 2^k} e^{-2^kt}
 \biggl(\int_{t_{k-1}}^te^{(s-t)2^{2k}}\,ds\biggr)\widehat w_k(\xi).
\end{eqnarray*}
In the second inequality we used our induction assumption.
Now, for all $t\ge t_k$, we have $t-t_{k-1}\ge 2^{-2k}\log2$,
so that $1-e^{(t_{k-1}-t)2^{2k}}\ge \frac{1}{2}$.
This in turn implies
$\int_{t_{k-1}}^te^{(s-t)2^{2k}}\,ds \ge 2^{-2k-1}$.
Then, for all $t\ge t_k$, we get
$$ \widehat u_t(\xi)\ge A^{2^{k}}  2^{2k+7-2^k} e^{-2^kt},$$
and~\eqref{induction} follows.

Moreover,
$\|\widehat w_k\|_{L^1}=(2\pi)^{-d} \|\widehat w_{k-1}\|_{L^1}^2$.
Since $\|\widehat w_0\|_{L^1}=1$, by induction we get
$$ \|\widehat w_k\|_{L^1}=(2\pi)^{-d(2^k-1)}.$$
Set $t^*=\lim_{k\to\infty} t_k=\log(2^{1/3})$.
We have $\widehat \psi(2^{-k}\,\cdot)\ge1$ in $E_k$.
Hence,
$$\widehat \psi(2^{-k}\xi)\widehat u_{t^*}(\xi)\ge 
  A^{2^k}\alpha_k(t^*)\widehat w_k(\xi)\ge0.$$  
Then
\begin{eqnarray*}
\|u_{t^*}\|_{\dot B^{a,\infty}_\infty} &\ge &
  \sup_{k\in\N} 2^{(a+d)k} |\psi(2^k\cdot)*u_{t^*}(0)| \\
&=&
  (2\pi)^{-d}\sup_{k\in\N} 2^{ak} \|\widehat \psi(2^{-k}\cdot)\widehat u_{t^*}\|_{L^1} \\
&\ge&
      \sup_{k\in\N} A^{2^k} (2\pi)^{-d\,2^k} 2^{ak+2k+7-2^k}e^{-2^kt^*}\\
&=& \sup_{k\in\N} \biggl\{\biggl(A\,2^{-4/3}(2\pi)^{-d} \biggr)^{2^k} 2^{(a+2)k+7}\biggr\}     
      =\infty.
\end{eqnarray*}
\endProof

\bigskip

\begin{remark}
\begin{rm}
The same proof goes through for solutions of~\eqref{general syst},
under the additional conditions $\widehat c_{j,h,k}\ge0$
for all $j,h,k=1,\ldots,m$, and, say, $\widehat c_{1,h,k}\ge c>0$.
In this case, one can obtain a solution $u=(u_1,\ldots,u_m)$,
such that the first component blows up, starting from
a datum $u_0$ such that $\widehat u_0\ge0$ and with the first component
satisfying the conditions of Theorem~\ref{theorem2}.
\end{rm}
\end{remark}

\begin{remark}
\begin{rm}
Analogous results can be obtained for space-periodic solutions of \eqref{general syst}, 
i.e. those defined on $d$-dimensional torus. Instead of the Fourier transform, we will consider the 
Fourier coefficients $\widehat u(\xi)$, $\xi\in \mathbb Z$. 
\end{rm}
\end{remark}

\bigskip

\noindent
{\bf Acknowledgments.}
This research was supported by the \'EGIDE--KBN POLONIUM project 6215.II/2005/2006, 
KBN (MNI) grant 2/P03A/002/24,  and 
by the European Commission Marie Curie Host Fellowship 
for the Transfer of Knowledge ``Harmonic Analysis, Nonlinear
Analysis and Probability''  MTKD-CT-2004-013389. 
\bibliographystyle{amsplain}

\end{document}